\documentclass[11pt]{article}
\usepackage{latexsym, amssymb, amsmath, amsthm, xspace,enumerate,color}

\theoremstyle{plain}
\newtheorem{theorem}{Theorem}[section]
\newtheorem{proposition}[theorem]{Proposition}
\newtheorem{lemma}[theorem]{Lemma}
\newtheorem{corollary}[theorem]{Corollary}
\theoremstyle{definition}

\newcommand{\GAP}{{\sf GAP}\xspace}
\newcommand{\MAGMA}{{\sc MAGMA}\xspace}
\newcommand{\F}[1]{{\mathbb F}_{#1}}
\newcommand{\nm}{\!-\!}

\def\Aut{{\rm Aut}}
\def\Hom{{\rm Hom}}
\def\Stab{{\rm Stab}}
\def\Out{{\rm Out}}

\def\C{{\sf C}}
\def\Sym{{\rm Sym}}

\def\Alt{{\rm Alt}}

\def\Alt{{\rm Alt}}
\def\GL{{\rm GL}}

\def\PGammaL{{\rm P \Gamma L}}

\def\PSL{{\rm PSL}}
\def\Sp{{\rm Sp}}

\def\GO{{\rm GO}}

\def\Sym{{\rm Sym}}

\def\Alt{{\rm Alt}}

\def\GL{{\rm GL}}

\def\F{{\mathbb F}}
\newcommand{\cC}{\mathcal{C}}
\newcommand{\cS}{\mathcal{S}}
\def\nonab{{\rm nonab}}
\def\sol{{\rm sol}}

\title{Comparing the order and the minimal number of generators of a finite
	irreducible linear group}
\author{Derek Holt and Gareth Tracey}

\begin{document}
	\maketitle
\begin{abstract}
We prove that $d(G) \log |G| = O(n^2 \log q)$ for irreducible subgroups $G$ of
$\GL(n,q)$, and estimate the associated constants. The result is motivated
by attempts to bound the complexity of computing the automorphism groups
of various classes of finite groups.
\end{abstract}
	\section{Introduction}\label{sec:intro}
	
	All groups in this paper are finite.
	For a group $G$, we denote the smallest
	size of a generating set of $G$ by $d(G)$.
	All logarithms are to the base $2$.

The complexity of a straightforward algorithm to compute the automorphism
group $\Aut(G)$ of a finite group $G$ is $O(|G|^{d(G)})$ and, frustratingly,
no tighter bound than this has been proved in general. This motivates
attempts to bound this quantity or, equivalently, its logarithm $d(G)\log |G|$
in the various classes of groups that arise in Computational Group Theory.
	
It is proved by the second author in \cite[Theorem 1.7]{Tracey18} that transitive subgroups $G \le \Sym(n)$
	satisfy $d(G)\log |G| = O(n^2/\sqrt{\log n})$.
	Here we prove an analogous result for finite irreducible linear groups.
	
	\begin{theorem}\label{thm:main}
		Let $G \le \GL(n,q)$ be irreducible. Then
\begin{enumerate}[\upshape 1.]
\item $d(G) \log |G| = O(n^2 \log q)$.
\item More precisely, we have $d(G)\log{|G|}\le Dn^2\log{q}$, where $D:=
C +\frac{10\log{3}}{9}+\frac{1}{2} \simeq 4.39$,
and $C = \frac{3\log{96}}{4\log{5}}$.
\item If, in addition, $G$ is (weakly quasi)-primitive, then Part 2
holds with $D = 2$, and this result is optimal. 
\end{enumerate}
	\end{theorem}
	
The proof is split into two parts in which the primitive and imprimitive cases
are handled.
In particular, we prove a stronger version of Part 3 in
Proposition~\ref{prop:primitive}.

For the benefit of readers who are only interested in
the bounds as orders of magnitudes, we shall prove them initially in this form,
in Sections~\ref{sec:primitive} and ~\ref{sec:general},
and then proceed to estimate the constants, which we do in
Section~\ref{sec:constants}. We suspect that the optimal bound
for $D$ is $2$. We could probably improve our current estimate by a more
detailed study of critical situations, but we see no immediate prospect
of getting close to the best possible result, so this does not seem worthwhile.

As a first step in bounding the constants, Section~\ref{sec:complen}
is devoted to the proof of Theorem~\ref{PrimBound},
which establishes a bound on the composition length of a
(weakly quasi)-primitive linear group over a finite field,
and may be of independent interest.

Note that the optimal asymptotic upper bounds for $d(G)$ and $\log |G|$
are respectively $O(n\log q/\sqrt{\log n})$ \cite[Theorem A]{LucMenMor00}
and $O(n^2 \log q)$ (the bound for $\GL(n,q)$). Our main result relies on
the fact that, in irreducible linear groups, when one of these two
terms is large, the other is relatively small.
	
We thank Colva Roney-Dougal for initially bringing this problem to our
attention. Indeed, part of our original motivation for Theorem~\ref{thm:main}
was to improve the bounds in \cite{RDSiccha} for the complexity of computing
normalizers in the symmetric group of affine-type primitive permutation groups
from $2^{O(\log^3 n)}$ to $2^{O(\log^2 n)}$ (though the authors of
\cite{RDSiccha} ultimately found a more direct and straightforward way of
achieving this).
We are grateful also to David Craven for assistance with the proof of
	Lemma~\ref{lem:tprep}.
	
	\section{The primitive case}\label{sec:primitive}

We start with some preliminary lemmas, and then deal with the weakly
quasiprimitive case of Theorem~\ref{thm:main} in
Proposition~\ref{prop:primitive}. The first lemma is attributed to Peter
Neumann.

\begin{lemma}\label{lem:dpg} \cite{CST89}
Let $G \le \Sym(t)$. Then $d(G) \le t/2$, except for $d(\Sym(3))=2$.
\end{lemma}

	\begin{lemma}\label{lem:tprep}
		Let G = $H \circ K$ be a central product of subgroups $H$ and $K$, and let $\rho$ be
		an absolutely irreducible representation of $G$ over a finite field $F$.
		Then $\rho$ is equivalent to $\rho_1 \otimes \rho_2$, where $\rho_1$ and
		$\rho_2$ are absolutely irreducible representations of $H$ and $K$,
		respectively.
	\end{lemma}
	\begin{proof} (This is also proved in \cite[3.16]{Aschbacher}.)
		This result  would follow from \cite[Chapter 3, Theorem 7.2]{Gorenstein} if the
		field $F$ were a splitting field for each of the central factors $H$ and $K$
		of $G$.
		
		If not, then let $E$ be a finite extension field of $F$ which does have this
		property. Then the natural extension $\rho^E$ of $\rho$ to a representation
		over $E$ is absolutely irreducible (since $\rho$ is absolutely irreducible),
		and so it is equivalent to $\rho_1 \otimes \rho_2$ for absolutely irreducible
		representations $\rho_1$ and $\rho_2$ of $H$ and $K$ over $E$.
		We need to prove that $\rho_1$ and $\rho_2$ are equivalent to representations
		over $F$, since the result will then follow from \cite[Theorem 29.7]{CurRei88},
		according to which two equivalent representations of a group over a field $E$
		that can both be written over a subfield $F$ of $E$ are equivalent as
		representations over $F$.
		
		By \cite[Theorem 74.9]{CurRei87}, an absolutely irreducible representation
		of a finite group over a field $E$ of characteristic $p>0$ can be written over
		a subfield $F$ if and only if the traces of the images of all group elements
		under the representation lie in $F$. Now the restriction $\rho_H$ of $\rho$
		to the subgroup $H$ is equivalent to $k\rho_1$, where $k$ is the degree of
		$\rho_2$, so $k\rho_1$ is equivalent to a representation over $F$.
		So the condition on traces is satisfied by $k\rho_1$ and hence also by
		$\rho_1$, and similarly for $\rho_2$. This completes the proof.
	\end{proof}
	
	Aschbacher \cite{Aschbacher} showed that maximal subgroups of
	classical groups over finite fields lie in one of nine classes, which
	he called $\cC_1$ -- $\cC_8$ and $\cS$.
	We extend (or abuse) this notation by applying it to arbitrary
	subgroups of maximal subgroups in classes $\cC_1$ -- $\cC_7$.
	For example, we
	view every semilinear subgroup of $\GL(n,q)$  a member of class $\cC_3$.
	Viewed in this way, the main result of \cite{Aschbacher} is that
	every subgroup $G$ of $\GL(n,q)$ either lies in at least one of
	$\cC_1$ -- $\cC_7$, or it lies in $\cC_8$ or $\cS$, in which case $G$
	has a normal absolutely irreducible subgroup
	that is simple modulo scalars, and $G$ does not lie in class $\cC_5$
	(that is, it is not defined over a proper subfield of $\F_q$ modulo scalars).
	
	For divisors $e$ of $n$, there are natural embeddings of the group
	$\Gamma L(n/e,q^e)$ into $\GL(d,q)$, and we shall denote their images
	by $\Gamma L(n/e,q^e)$. (These images are all conjugate in $\GL(d,q)$.)
	The semilinear subgroups of $\GL(n,q)$ are exactly the conjugates of subgroups
	of $\Gamma L(n/e,q^e)$ for some $e>1$.
	
	The following lemma describes a useful condition for a group $G$ to be
	semilinear. This is proved in a much more general setting in
	\cite[11.5]{Aschbacher}. There is also an informal proof in
	\cite[Section 2]{Smash}, which describes an algorithm for computing
	the associated embedding of $\Gamma L(n/e,q^e)$ into $\GL(d,q)$.
	We shall give a sketch proof here, which is essentially the same as the
	one in \cite{Aschbacher}.
	
	\begin{lemma}\label{lem:c3sub}
		Let $N \unlhd G$, let $\phi:G \to \GL(n,q)$ be a representation,
		and suppose that the associated $\F_qN$-module $V$ is homogeneous, and that
		the centralizing field of its  irreducible constituents is $\F_{q^e}$.
		Then $\phi(N) \le \GL(n/e,q^e)$ with $\phi(G) \le \Gamma L(n/e,q^e)$ (so
		$G$ is semilinear).
		
		Furthermore, the action of $\phi(N)$ on its natural module is homogeneous
		with absolutely irreducible constituents.
	\end{lemma}
	\begin{proof}
		Let $C := C_{\GL(n,q)}(\phi(N))$, $H := \langle C,\phi(N) \rangle$, and
		$Z := Z(C) = Z(H)$.
		By \cite[Chapter 3, Theorem 5.4\,(iii)]{Gorenstein} (the assumption that
		$\Hom_G(V_1,V_1)=F$ is not needed in the proof) or \cite[3.10\,(2)]{Aschbacher},
		we have $C \cong \GL(t,q^e)$, where $t$ is the number of irreducible
		consituents of $V$. Furthermore, by \cite[Chapter 3, Theorem 5.5]{Gorenstein}
		(which also does not make use of the assumption $\Hom_G(V_1,V_1)=F$),
		$C$ acts transitively on the set of irreducible $\F_qN$-submodules of $V$,
		so $H$ acts irreducibly, and its centre $Z$, which is cyclic of order $q^e-1$,
		must act homogeneously. 
		
		Since $\F_{q^e}$ is also the centralizing field of the constituents of the
		associated $\F_qZ$-module, we have $C_{\GL(n,q)}(Z) \cong \GL(n/e,q^e)$ and,
		since $C_{\GL(n,q)}(Z)$ acts irreducibly, it can be identified with the image
		of the natural embedding of $\GL(n/e,q^e)$ into $\GL(n,q)$.
		So, in particular, we have $\phi(N) \le \GL(n/e,q^e)$.
		Since $\phi(G)$ normalizes $\phi(N)$, it must also normalize $C$, $Z$, and
		$C_{\GL(n,q)}(Z)$, so $\phi(G) \le \Gamma L(n/e,q^e)$ as claimed.
		
		The fact that $\phi(N)$ acts homogeneously with absolutely irreducible
		constituents as a subgroup of $\GL(n/e,q^e)$ follows easily from the fact
		that its centralizer in $\GL(n/e,q^e)$ is isomorphic to $\GL(t,q^e)$.
	\end{proof}

	We recall from~\cite{HRD13} that, for a field $K$, a subgroup of $\GL(n,K)$
	is said to be {\it weakly quasiprimitive} if all of its characteristic
	subgroups act homogeneously. A{\it subdirect subgroup} of a direct product
	of groups is defined to be a subgroup that projects onto each of the direct
	factors.
	
	\begin{lemma}\label{lem:slred}
		Let $G \le \GL(n,q)$ be weakly quasiprimitive. Then there
		exists a divisor $e \ge 1$ of $n$ such that $G \le \Gamma L(n/e,p^e)$,
		and $K := G \cap \GL(n/e,q^e)$ has the properties listed below.
		We regard $K$ as a subgroup of $\GL(n/e,p^e)$, and denote its associated module
		of dimension $n/e$ over $\F_{q^e}$ by $V$.
		
		\begin{enumerate}
			\item $|G:K|$ divides $e$, and $K$ is either normal or 2-step subnormal
			in $G$, with all associated quotient groups cyclic.
			\item $Z(K) = K\cap Z(\GL(n/e,q^e))$, and $K/Z(K)$ is isomorphic to a subdirect
			subgroup of $K_1\times \cdots\times K_u$, where each $K_i$ satisfies
			one of the following:
			\begin{enumerate}[\upshape(a)]
				\item $K_i$ has shape $N.S$, where $N$ is elementary abelian of order
				$r_i^{2m_i}$, with $r_i$ a prime dividing  $q^e-1$, and $S$ is a
				completely reducible subgroup of $\Sp(2m_i,r_i)$ acting naturally
				on $N$; or
				\item there exists $s_i \ge 2$ and a finite simple group $S$ with an absolutely irreducible
				projective representation $\rho$ of dimension $s_i$ over
				$\mathbb{F}_{q^e}$,
				and a subgroup $A$ of the stabilizer of $\rho$ in $\Aut S$, such that
				$K_i$ is a subgroup of $A\wr \Sym(t_i)$ containing $S^{t_i}$.
			\end{enumerate}
			Moreover, setting $f(K_i):=r_i^{m_i}$ if $K_i$ is as in (a),
			and $f(K_i):=s_i^{t_i}$ if $K_i$ is as in (b), we have that
			$\prod_{i=1}^uf(K_i)$ divides $n/e$. Finally, if $r$ is a fixed prime,
			then there is at most one value of $i$ in the range $1\le i\le u$ with the
			property that $K_i$ is as in (a) with $r_i=r$. 
		\end{enumerate}
	\end{lemma}
	\begin{proof}
		The proof is an adaptation of that of \cite[Theorem 1.2, Part 2]{HRD13}.
		Essentially the same results are established in
		\cite[Lemmas 5.1 and 5.2, and the proof of Corollary 5.3]{Tracey21}
		
		By \cite[proof of Lemma 3.1]{LucMenMor01}, if $G$ has an abelian non-scalar
		characteristic subgroup, then there is a divisor $d>1$ of $n$ such that $G$
		has a characteristic subgroup $M$ that can be identified with a subgroup
		of $\GL(n/d,q^d)$, where $G/M$ is cyclic of order dividing $d$.
		(So $G$ is semilinear.) Furthermore, $M$ is weakly quasiprimitive and has no
		abelian nonscalar characteristic subgroups.  In particular, the centre
		$Z:=Z(M)$ of $M$ is equal to the scalar subgroup of $M$.
		
		Let $L$ be the generalized Fitting subgroup $F^*(M)$ of $M$. Then $L$ is a
		central product of the non-central subgroups $O_{r_i}(M)$ for a (possibly empty)
		list of distinct primes $r_1,r_2,\ldots,r_v$,
		and a (possibly empty) list of normal subgroups
		$T_{v+1},T_{v+2},\ldots,T_u$ of $M$, where each $T_i$ is
		a central product of $t_i$ copies of a quasisimple group $S_i$, and $M$
		(acting by conjugation) permutes these copies transitively.
		By the standard properties of $F^*(M)$, we have $C_M(F^*(L)) = Z(F^*(L))$ and,
		since this is a characteristic abelian subgroup of $M$, it must be equal
		to $Z$.
		
		Since the groups $O_r(M)$ with $r=r_i$  are non-central and have no non-central
		abelian characteristic subgroups, it follows from \cite[Lemma 1.7]{LucMenMor01}
		that $O_r(M)$ is the central product of a (cyclic) Sylow $r$-subgroup of $Z$
		and an extraspecial group with centre contained in $Z$. Furthermore,
		we have $$C_M(O_r(M/Z))/C_M(O_r(M))= O_r(M)Z/Z,$$
		and $M/(O_r(M)C_M(O_r(M))$  is isomorphic to a completely reducible subgroup
		of $\Sp(2m,r)$, where $|O_r(M)Z/Z| = r^{2m}$.
		
		Since $M$ is weakly quasiprimitive, $L$ is homogeneous. If its restriction
		to its irreducible constituents is not absolutely irreducible then, by
		by Lemma~\ref{lem:tprep} there is a divisor $e/d$ of $n/d$, and a normal
		subgroup $K \le \GL(n/e,q^e)$ of $M$ with $L \le K$ and $M/K$ cyclic of order
		dividing $d/e$. So we have proved that $K$ satisfies Condition 1 in
		the statement of the lemma.
		
		Note that, since $L \le K \unlhd M$, for each $i \le v$ we have
		$O_{r_i}(K)=O_{r_i}(M)$ and $K/(O_r(K)C_K(O_r(K))$ is isomorphic to a
		completely reducible subgroup of $\Sp(2m,r)$. It is possible that $K$ fails to
		act transitively on the  $t_i$ copies of $S_i$ for some $i>v$, and in that case
		we increase $u$ and redefine the subgroups $T_i$ to reimpose this condition.
		Note also that $Z(K) \le C_M(L)$, which is contained in the scalar subgroup
		of $\GL(n/d,q^d)$, so $Z(K) = K \cap Z(\GL(n/e,q^e))$.
		
		By Lemma~\ref{lem:tprep}, $L$ remains homogeneous as a subgroup of
		$\GL(n/e,q^e)$ and acts absolutely irreducibly on its irreducible constituents.
		Let $V$ be one such constituent. So $\dim V$ divides $n/e$.
		It follows from Lemma~\ref{lem:tprep} that $V$ is isomorphic to a tensor
		product of a one-dimensional $Z(K)$-module, absolutely irreducible modules
		$V_i$ for each $O_{r_i}(K)$, and absolutely irreducible modules $V_{T_i}$
		for each $T_i$.
		
		Note that, if the lists of central factors $O_{r_i}(K)$ and $T_i$ of $L$  are
		both empty, then $n=e$ and there is nothing to prove. So we may assume that
		at least one of them is non-empty.
		
		We define $K_i =  K/C_K(O_{r_i})$ for $1 \le 1 \le v$ and
		$K_i = K/C_K(T_i)$ for $v+1 \le i \le u$. Since $C_K(L) = Z(L) =
		K \cap Z(\GL(n/e,q^e)$, the natural map $G \to K_1 \times \cdots \times K_t$
		has kernel $Z(K)$ and its image is a subdirect subgroup of
		$K_1 \times \cdots \times K_u$.
		
		We have already proved that $K_i$ has the shape described in Condition 2\,(a)
		in the statement of the lemma for $1 \le i \le u$. 
		By \cite[Chapter 5, Theorem 5.5]{Gorenstein}, the faithful absolutely
		irreducible representations of $O_{r_i}(K)$ have dimension $r_i^{m_i}$,
		and we define $f(K_i) = r_i^{m_i}$. So $\dim V_i = f(K_i)$.
		Since $Z(O_{r_i}(K))$ has order $r_i$ and is a subgroup of the scalar
		subgroup $Z(K)$, we have $r_i|q^e-1$.
		
		Lemma~\ref{lem:tprep} tells us that, for $i>u$, the module $V_{T_i}$ 
		decomposes as a tensor product of absolutely irreducible representations
		of its quasisimple central factors. Since these central factors are permuted
		transitively under the conjugaction action of $L$, they must all have the
		same dimension $s_i$ with $s_i \ge 2$, and then $\dim V_{T_i} = s_i^{t_i}$,
		and we define $f(K_i) = s_i^{t_i}$.
		
		So now we have $\dim V = \prod_{i=1}^f(K_i)$, and this divides $n/e$,
		which completes the proof of the lemma.
	\end{proof}
	
	\begin{proposition}\label{prop:primitive}
		Let $H$ be a subnormal subgroup of a weakly quasiprimitive group
		$G \le \GL(n,q)$.  Then $d(H) \log |G| \le Cn^2 \log q$ for some constant $C$.
		
			In fact this result is true with $C = \frac{3 \log(96)}{4\log(5)} \simeq 2.127$,
			and $d(G) \log |G| \le 2 n^2 \log q$. Both of these results are optimal.
	\end{proposition}
	\begin{proof}
		The bound with $C = \frac{3 \log(96)}{4\log(5)}$ is
		attained with $G$ equal to the normalizer of a $2$-group $H$ of symplectic
		type (with $d(H)=3$) in $\GL(2,5)$. The optimality of
		$d(G) \log |G| \le 2n^2 \log q$  arises from the groups $G=\GL(n,q)$ for
		large $n$.
		
As we said earlier, we shall initally prove only that the orders of magnitude
are as claimed, and postpone the estimation of the constants until
Section~\ref{sec:constants}.  
		
		Suppose that $H$ is a subnormal subgroup of the weakly quasiprimitive group
		$G \le \GL(n,q)$. Since $G$ acts homogeneously, we may assume that it acts
		irreducibly. The result is clear if $n=1$, so assume that $n>1$.
		
		We apply Lemma~\ref{lem:slred} to the group $G$.
		If the result holds for the subnormal subgroup $H \cap K$ of $K$, then
		\begin{eqnarray*}
			d(H) \log |G| &\le &(d(H \cap K) + 2)(\log |K| + \log e)\\
			&\le& C\frac{n^2}{e} \log q + (d(H \cap K)+2)\log e + 2\log |K|.
		\end{eqnarray*}
		Now $d(H \cap K) \le 2\log(n/e) + 1$ by \cite[Theorem 1.2, Part 2]{HRD13},
		and $\log |K| \le \frac{n^2}{e} \log q$, so
		$$d(H) \log |G| \le \frac{n^2(C+2)}{e}\log q +
 \left(2\log \frac{n}{e}+3 \right) \log e = O(n^2 \log q). \eqno{(1)}$$

		So we can assume that $K=G$ and $e=1$ in the conclusion of
		Lemma~\ref{lem:slred}, and we shall rename the direct factors
		$K_i$ of $G/Z(G)$ as $G_i$.  We shall refer to the two possibilities for
		$G_i$ that are described in Cases 2\,(a) and 2\,(b) of the conclusion
		of Lemma~\ref{lem:slred} as the solvable and non-solvable cases, respectively.
		
		Since $G/Z(G)$ is a subdirect subgroup of $G_1 \times \cdots \times G_u$ and
		$H$ is a subnormal subgroup of $G$, there are subnormal subgroups $H_i$
		of $G_i$ such that $d(H) \le  \sum_{i=1}^u d(H_i) +1$. We also have 
		$\log |G| \le \sum_{i=1}^u \log |G_i|+ \log (q \nm 1)$.
		In the proof of \cite[Theorem 1.2, Part 2]{HRD13}, it is proved that
		$d(H_i) \le 2 \log  f(G_i)$ in both the solvable and non-solvable cases.
		
		We shall prove first that
		\[
		(d(H_i)+1)(\log(q \nm 1) + \log |G_i|) \le C_1f(G_i)^2 \log q\ 
		\ \hbox{or}\  \  C_2f(G_i)^2 \log q \eqno{(2)}
		\]
		for suitable constants $C_1$ and $C_2$ in the solvable and non-solvable cases.

		Consider first the solvable case and, with the notation of 2\,(a) of
		Lemma~\ref{lem:slred}, put $r=r_i$ and $m=m_i$.
		Then $G_i = N.S$ with $N$ an elementary abelian $r$-group of order $r^{2m}$
		and $S \le \Sp(2m,r)$, where $r|q-1$, and that $f(G_i) = r^m$.
		Since $|\Sp(2m,r)| \le r^{2m^2+m}$, we have
		$|G_i| \le r^{2m}r^{2m^2+m}$, so
$\log |G_i| \le (2m^2+3m) \log r \le 2\left(\log f(G_i)\right)^2 + 3\log f(G_i)$. 
		Then, since (as observed above) $d(H_i) \le 2 \log f(G_i)$, we have
\begin{eqnarray*}
	(d(H_i) + 1)(\log(q \nm 1)+\log |G_i|) & \le &\\
(2 \log f(G_i)+1)(\log (q \nm 1) + 2\left(\log f(G_i)\right)^2 + 3\log f(G_i)))
	&\le &\\ C_1 f(G_i)^2 \log q && \quad\quad\quad\quad(3)
\end{eqnarray*}
		for some constant $C_1$.

		Consider now the non-solvable case and, with the notation of 2\,(b) of
		Lemma~\ref{lem:slred}, put $s=s_i$ and $t=t_i$.
		Then $S^{t} \le G_i \le A \wr \Sym(t)$ where $A= \Aut S$ with $S$
		nonabelian simple, and $f(G_i)=s^{t}$, where $S_i$ has an absolutely
		irreducible projective representation of degree $s$ over $\F_{q}$.
	So $|G_i|\leq q^{ts^2}t!$ and $\log|G_i| \le ts^2 \log q + t\log t$.
		
		Now, since $d(H_i) \le 2 \log f(G_i)$, when $t >1$ we have
		\begin{eqnarray*}
			(d(H_i)+1)(\log(q \nm 1)+\log |G_i|) &\le &\\
			(2t \log s + 1) ((1+ts^2)\log q + t\log t)
			&\le & C_2 s^{2t} \log q = C_2 f(G_i)^2 \log q
                        \quad(4)
		\end{eqnarray*}
		for some constant $C_2$.

		When $t=1$, the group $G_i$ is almost simple and so $H_i$ is either
		trivial or almost simple, and we have $d(H_i) \le 3$.  So
$$(d(H_i)+1)(\log (q \nm 1) + \log |G_i|) \le 4f(G_i)^2 \log q. \eqno{(5)}$$

	This completes the proof of our claims in (2) for the individual factors $G_i$
		and of the theorem statement in the case $u=1$,
	and we now proceed to the proof of the theorem in the case $u>1$.
		
	So we have proved already that
	$(d(H_i) + 1)(\log(q \nm 1) + \log |G_i|) \le \max(C_1,C_2) f(G_i)^2 \log q$
		for $1 \le i \le u$.
		
	Now, using $d(H_i) \le \lfloor 2\log f(G_i)  \rfloor \le f(G_i)$ and
$\log |G_i| \le f(G_i)^2 \log q$, we have
		\begin{eqnarray*}
			d(H) \log|G|& \le &
\left(1 + \sum_{i=1}^{u} d(H_i)\right)\left(\log(q\nm 1)+
                        \sum_{i=1}^{u} \log |G_i|\right)\\ & \le &
   \sum_{i=1}^{u} \left(1 + d(H_i)\right)\left(\log(q\nm 1)+\log |G_i|\right)
+ \sum_{\underset{i \ne j}{i,j=1}}^{u} d(H_i)\log |G_j|\\
& \le & \left(\max(C_1,C_2)\sum_{i=1}^{u} f(G_i)^2 +
\sum_{\underset{i \ne j}{i,j=1}}^{u} f(G_i) f(G_j)^2\right)\log q.
		\end{eqnarray*}

Now, since $\prod_{i=1}^u f(G_i)$ divides $n$ and $f(G_i) \ge 2$ for all $i$,
we have $f(G_i) \le n/2^{u-1}$ for all $i$, and
$f(G_i) f(G_j)^2 \le n^2/2^{2u-3}$ for all $i \ne j$, and so
\[
d(H) \log|G|  \le  \left (\frac{\max(C_1,C_2)u}{2^{2u-2}} +
\frac{u(u-1)}{2^{2u-3}}\right)n^2\log q, \eqno{(6)}
\]
which is at most $Cn^2\log q$ for sufficiently large $C$.

	\end{proof}
	
\section{The general case}\label{sec:general}
In this section we complete the proof of Part 1 of Theorem~\ref{thm:main}.
A subgroup of a wreath product $H \wr T$, with $T$ a permutation group, is
said to be {\em large} if it projects onto $T$ and if the stabiliser of a
block in the natural imprimitive action projects onto $H$.
For a group $G$, we denote the length of a composition series of $G$ by $a(G)$.

Let $G\le \GL(n,q)$ be irreducible.
It is proved in \cite{Suprunenko76} that $G$ can be identified
with a large subgroup of $H \wr T$ with $H \le \GL(s,q)$ primitive and
$T \le \Sym(t)$ transitive, where $n=st$ and either $H = G$ is primitive and
$t = 1$; or $G$ is imprimitive with minimal blocks of dimension  $s$, $H$
is the action of a block stabiliser on that block, and $T$ is the induced
action of $G$ on the set of blocks.

When $t=1$, the group $G$ is primitive, and the result follows from
Proposition~\ref{prop:primitive}. So assume now that $t>1$.

Let $K \le H^t$ be the kernel of the action of $G$ on the blocks of the
wreath product decomposition.
Let $H_1$ be the projection of $K$ onto the first of the direct factors
of $H^t$ and, for $2 \le i \le t$, let $H_i$ be the projection onto the
$i$th direct factor of the subgroup of $K$ that induces the identity on each
of the first $i-1$ factors. So each $H_i$ is isomorphic to a subnormal
subgroup of the primitive group $H \le \GL(s,q)$.

We have $G/K \cong T$, and so $\log |G| = \log |K| + \log |T|$,
and hence $d(G) \log |G| \le d(G) \log |K| + d(G) \log |T|$.
We shall prove the required bound for the two terms
$d(G) \log |K|$ and $d(G) \log |T|$ separately.

We start by proving that $d(G) \log |T| = O(n^2)\log q$, for which the proof is
similar to that of \cite[Theorem 1.7]{Tracey18}.

Now $T$ can be identified with a subgroup of $T_1 \wr  \cdots \wr T_m$,
where each $T_i$ is a primitive permutation group of degree $t_i > 1$, and
$\prod_{i=1}^m t_i = t$. If none of the the groups $T_i$ with $t_i \ge 3$
contains $\Alt(t_i)$ then, by \cite[Corollary 1.4]{Maroti02}, we have
$|T| \le c^t$ with $c = 2.8349$.
Furthermore, since by \cite[Theorem 1.2]{HRD13}, we have
$d(H_i) \le 2 \log s + 1$ for each of the subnormal subgroups $H_i$ of $H$
and $d(T) < t$ by Lemma~\ref{lem:dpg},
we have $d(G) \le t(2 \log s + 1) + d(T) < 2t(\log s + 1)$,
and so $d(G) \log |T|  = O(t^2\log s) = O(n^2)$, as claimed.

Otherwise, choose $k$ with $1 \le k \le m$ such that $t_k$ is maximal with
$t_k \ge 3$ and $\Alt(t_k) \le T_k$. Then,
by \cite[Proposition 6.9(i)]{Tracey18}, we have $d(G) < a(R)v$, where
$R = H \wr T_1 \wr \cdots T_{k-1}$, $S = T_k \wr \cdots \wr T_m$, and
$v = \prod_{i=k+1}^m t_i$.

Also, by \cite[Theorem C]{LucMenMor01}, we have
\[a(R) = O(\deg R\log q) =
O(s\prod_{i=1}^{k-1} t_i \log q) = O\left(\frac{n \log q}{t_kv}\right),\]
and so $d(G) = O\left(\frac{n \log q}{t_k}\right)$.
Now, by \cite[Proposition 1.9]{Tracey16}, we have $|T| \le t_k^t$, and
then $d(G) \log |T| = O(\frac{nt}{t_k} \log t_k\log q) =
O(n^2\log q)$ as claimed.

Finally we prove that $d(G) \log |K| = O(n^2 \log q)$.  Now $d(G) \le
d(K) + d(T) \le \sum_{i=1}^s d(H_i) + (t+1)/2$ by Lemma~\ref{lem:dpg}, and
$\log |K| \le t \log |H|$ so by Proposition~\ref{prop:primitive} we have
\begin{align*}
d(G) \log |K| &\le \left(\sum_{i=1}^t d(H_i) + \frac{t+1}{2}\right)t\log |H|\\ 
&\le Cs^2t^2\log q + s^2t\frac{t+1}{2} \log q = O(n^2 \log q),&\quad(7)
\end{align*}
which completes the proof.

\section{Composition length in weakly quasiprimitive linear groups}\label{sec:complen}
	Recall that for a finite group $G$, we write $a(G)$ for the composition length of $G$. As mentioned in the introduction, the purpose of this section is to find upper bounds on $a(G)$ when $G$ is a weakly quasiprimitive subgroup of $\GL(n,q)$. To state our main result, we note that the group $\GL(4,2)$ has a weakly
quasiprimitive subgroup $X$ of shape $\Sym(3)\wr\Sym(2)$ with $a(X)=5$, and that
$X$ has three subgroups $X_1,X_2$ and $X_3$ of index $2$ which are also weakly quasiprimitive. Our main result can now be stated as follows.
	\begin{theorem}\label{PrimBound} Let $q$ be a power of a prime, and let $G\le \GL(n,q)$ be weakly quasiprimitive. Then $a(G)\le \log{n}+n\log{q}$,
		except when either: 
		\begin{enumerate}[\upshape(i)]
\item $(n,q)=(2,3)$ and $G=\GL(2,3)$, in which case $a(G)=5$; 
\item $(n,q)=(2,5)$ and $G$ has shape $4\circ 2^{1+2}.\Sym(3)$, in which case $a(G)=6$; or
\item $(n,q)=(4,3)$ and $G$ has shape $G=2^{1+4}.L$, where either $L\cong X$
(as defined above), in which case $a(G)=10$;
or $L\in\{X_1,X_2,X_3\}$, in which case $a(G)=9$.
		\end{enumerate}
	\end{theorem}
	
	Note that the bound in Theorem \ref{PrimBound} is almost best possible. Indeed, if $n$ is a power of $2$, then the semilinear group $G:=\Gamma \mathrm{L}(1,2^n)\le \GL(n,2)$ is weakly quasiprimitive of shape $(2^n-1):n$, and composition length $\omega(2^n-1)+\log{n}$ where $\omega(m)$ is the number of prime divisors of $m$, counted with multiplicities.
	
To prove Theorem \ref{PrimBound}, we shall start by proving the following two
propositions, and then use these together with Lemma~\ref{lem:slred} to
complete the proof.
%
	\begin{proposition}\label{prop:rmCase}
		Let $G$ be a group of shape $N.L$, where $N$ is an elementary abelian group of order $r^{2m}$, with $r$ prime, and $L$ is a completely reducible subgroup of $\Sp(2m,r)$ acting naturally on $N$. Suppose that Theorem \ref{PrimBound} holds in dimensions $n\le 2m$, and let $q$ be a prime power with $r$ dividing $q-1$.
Then we have:
	\begin{enumerate}[\upshape(i)]
\item If $(m,r,q)\neq (1,2,3),(1,2,5),(1,2,7),(1,3,4),(2,2,3)$, or $(2,2,5)$, then $a(G)\le m\log{r}+(r^m-1)\log{q}$.
In the exceptional cases, we have $a(G)=4,4,4,6,9$, or $9$, respectively.
\item If $(m,r,q)\not\in\{(1,2,3),(2,2,3)\}$,
                 then $a(G)\le m\log{r}+r^m\log{q}$.
		\end{enumerate}
	\end{proposition} 
	
	\begin{proposition}\label{prop:stCase}
		Let $S$ be a nonabelian finite simple group, let $q$ be a prime power, and let $s$ be the dimension of a non-trivial projective irreducible representation $\rho$ of $S$ over $\mathbb{F}_q$. Also, let $A$ be a subgroup of $\Aut S$ which stabilizes $\rho$, and let $G$ be a subgroup of $A\wr \Sym(t)\cong A^t\rtimes\Sym(t)$ containing $S^t$. Then $a(G)\le t\log{s}+(s^t-1)\log{q}$.
	\end{proposition}

We begin our preparations toward the proofs of Propositions \ref{prop:rmCase}
and \ref{prop:stCase} with
an easy lemma concerning certain subgroups of direct
products of groups. Indeed, let $G_1,\hdots,G_t$ be groups, and write
$\pi_i:G_1\times\cdots\times G_t\rightarrow G_i$ for the natural
projection maps. Recall that a subgroup $G$ of the direct product
$G_1\times\cdots\times G_t$ is called a \emph{subdirect subgroup}
if $G\pi_i=G_i$ for all $i$.
	\begin{lemma}\label{lem:subdlemma}
		Let $G_1,\hdots, G_t$ be finite groups, and let $G$ be a subdirect subgroup of
		the direct product $G_1\times\cdots\times G_t$.
		Then $a(G)\le \sum_{i=1}^t a(G_i)$.
	\end{lemma}
	\begin{proof}
		We prove the lemma by induction on $t$. If $t=1$ then the result is trivial,
		so assume that $t\geq 2$. Let $K$ be the kernel of the projection map $\pi_t$.
		Then $K$ is isomorphic to a subdirect subgroup of the direct product
		$(K\pi_1)\times\cdots\times (K\pi_{t-1})$. The inductive hypothesis then implies that $a(K)\le \sum_{i=1}^{t-1}a(K\pi_i)$. Also, $K\pi_i\unlhd G\pi_i=G_i$. Thus, since $a(N)\le a(N)+a(X/N)=a(X)$ for finite groups $N$ and $X$ with $N\unlhd X$, we have $a(K)\le \sum_{i=1}^{t-1}a(G_i)$. Since $G/K\cong G_t$, the result follows.
	\end{proof}

	Next, we need an upper bound on composition lengths of subgroups of $\Sym(n)$, proved in \cite[Theorem 1.2]{GPRV}.
	\begin{proposition}\cite[Theorem 1.2]{GPRV}\label{prop:TransCompLength}
		Let $G\le\Sym(n)$ be a permutation group of degree $n$. Then $a(G)\le 4(n-1)/3$.
	\end{proposition}

	We remark that the bound in Proposition \ref{prop:TransCompLength} is best possible, even if one restricts to transitive groups. Indeed, as shown in \cite{GPRV}, if $n=4^b$, then the iterated wreath product of $b$ copies of $\Sym(4)$ embeds as a (transitive) subgroup of $\Sym(n)$, and has composition length $4(n-1)/3$. 
	
	The following corollary of Proposition \ref{prop:TransCompLength} will be needed for our proof of Proposition \ref{prop:rmCase}.
	\begin{corollary}\label{cor:PrimBoundcorollary}
		Let $q$ be a power of a prime $p$, and let $G\le \GL(k,q)$ be completely reducible.
Suppose that Theorem \ref{PrimBound} holds in dimensions less than or equal to $k$. Then 
\begin{enumerate}[\upshape(i)]
    \item $a(G)\le 4(k-1)/3+k\log{q}$ if $q\not\in\{3,5\}$;
    \item $a(G)\le 19k/6-4/3$ if $q=3$; and
    \item $a(G)\le 22k/6-4/3$ if $q=5$.
\end{enumerate}
In particular, $a(G)\le (7/3)k\log{q}$.
	\end{corollary}
	\begin{proof} Write $V$ for the natural $\mathbb{F}_q[G]$-module. Since $G$ is completely reducible, $V$ may be written in the form $V=V_1\oplus\cdots\oplus V_r$, where each $V_i$ is an irreducible $\mathbb{F}_q[G]$-module. It follows that $G$ embeds as a subdirect subgroup of $G^{V_1}\times\cdots\times G^{V_r}$, where $G^{V_i}$ is the irreducible linear group induced by $G$ on $V_i$. We deduce from Lemma \ref{lem:subdlemma} that $a(G)\le \sum_{i=1}^ra(G^{V_i})$. Since $k=\sum_{i=1}^r\dim{V_i}$, we may therefore assume that $G$ is irreducible. 
		
		So assume that $G$ is irreducible. The result follows quickly from the hypothesis of the corollary and Theorem \ref{PrimBound} if $G$ is primitive, so assume that $G$ is imprimitive, and let $W\le V$ be a minimal non-trivial block for $G$. Set $S:=\Stab_G(W)^W$ be the linear group induced by the (setwise) stabilizer $\Stab_G(W)$ of $W$ on $W$, and let $\Sigma$ be the set of $G$-translates of $W$. Write $m=\dim{W}$ and $t=|\Sigma|$, so that $k=mt$. Then $S\le \GL(m,q)$ is primitive, since $W$ is a minimal block for $G$, and $G^{\Sigma}\le \Sym(t)$ is transitive, since $G$ is irreducible. Moreover, $G$ embeds as a subgroup in the wreath product $S\wr G^{\Sigma}$ in such a way that the intersection $K:=G\cap S^{t}$ of $G$ with the base group of the wreath product is a subdirect subgroup in a direct product of $t$ copies of a normal subgroup $A$ of $S$. It follows from Lemma \ref{lem:subdlemma} that $a(K)\le a(A)t$,
and so $a(G)\le a(A)t+a(G/K)$.
Now Proposition \ref{prop:TransCompLength} applied to $a(G/K)$ give
 $a(G)\le a(A)t+4(t-1)/3$, and we can apply
the hypothesis of the corollary and Theorem \ref{PrimBound} to bound $a(A)$.

When $q\not\in\{3,5\}$, we find using some straightforward calculus that,
for fixed $q$ and $k$, and subject to the condition $mt=k$,
the function $t(\log{m}+m\log{q})+4(t-1)/3$ is
is decreasing in $m$ and is therefore maximized at $m=1$.
So $a(G)\le k \log q + 4(k-1)/3$, which completes the proof in this case.

When $q=3$ we find that
$a(G) \le t(\log{m}+m\log{q})+4(t-1)/3$ when $m\not\in\{2,4\}$;
that $a(G)\le 5t+4(t-1)/3=19k/6-4/3$ when $m=2$;
and that $a(G)\le 10t+4(t-1)/3=34k/12-4/3$ when $m=4$.
As in the previous paragraph above, it is easy to show that the largest of
these three quantities is always $19k/6-4/3$. The argument for $q=5$ is almost identical, and this completes the proof.
	\end{proof}
	
	We can see from the last paragraph of the proof of Corollary \ref{cor:PrimBoundcorollary}, and the remark after the statement of Proposition \ref{prop:TransCompLength}, that the bound in the corollary is best possible for $q=3$. Note also that the bound in Corollary \ref{cor:PrimBoundcorollary} is almost best possible for $q\neq 3$. Indeed, if $q=2^a+1$ is a Fermat prime, $k=4^b$ is a power of $4$, $Y$ is the iterated wreath product of $b$ copies of $\Sym(4)$, and
$X:=\GL(1,q)\cong \C_{q-1}$, then $G:=X\wr Y$ is an irreducible subgroup of $\GL(k,q)$ of composition length $ka+4(k-1)/3=k\log{(q-1)}+4(k-1)/3$.
	
	We can now prove Propositions \ref{prop:rmCase} and \ref{prop:stCase}.
	\begin{proof}[Proof of Proposition \ref{prop:rmCase}] Let $G$, $N$, $L$, $r$, $m$, and $q$ be as in the statement of the proposition.
		Then $a(G)=a(N)+a(L)=2m+a(L)$. If $(m,r)=(1,2)$, $(1,3)$, or $(2,2)$, then we can quickly see that, respectively, $a(L)\le 2,4$ or $5$, so $a(G)\le 4,6$ or $9$.
Parts (i) and (ii) follow in these cases, and it remains only to prove (i) in
the non-exceptional cases.
		
	Assume now that $(m,r)\not\in\{(1,2),(1,3),(2,2)\}$.
Suppose also that $r\neq 3$.
Then $a(G)=2m+a(L)\le 2m+8m/3+2m\log{r}$, by Corollary \ref{cor:PrimBoundcorollary}. Since $q\geq r+1$, this gives $a(G)\le \log{r^m}+(r^m-1)\log{q}$
for $(m,r,q)\not\in\{(3,2,3),(3,2,5)\}$. We can check by direct computation
in the completely reducible subgroups of $\Sp(6,2)$
that $a(L)\le 8$ for $(m,r)=(3,2)$. This gives $a(G)\le \log{r^m}+(r^m-1)\log{q}$ for $(m,r)\in\{(3,2,3),(3,2,5)\}$, and completes the proof of (i) for $r\neq 3$.
		
		Continuing to assume that $(m,r)\neq (1,3)$,
suppose now that $r=3$. Then $a(G)\le 2m+38m/6-4/3$, by Corollary \ref{cor:PrimBoundcorollary}. Since $r=3$ divides $q-1$, it is easy to see that this is less than $m\log{3}+(3^m-1)\log{q}$ whenever $m\geq 2$.
This completes the proof.
	\end{proof}
	
	\begin{proof}[Proof of Proposition \ref{prop:stCase}]
		Let $S$, $A$, $s$, and $t$ be as in the statement of the proposition. Then since $S$ embeds as a subgroup of $\mathrm{PGL}_s(q)$, $S$ has a non-trivial permutation representation of degree $(q^{s}-1)/(q-1)< q^s$. It then follows from \cite[Lemma 2.6]{JZPyb} that $|A/S|\le |\Out S|\le 3s\log{q}$. 
		
		Now, $G$ has a normal series $1<N<K<G$ where $N\cong S^t$, $K/N$ is a subgroup of $A^t$, and $G/K\le \Sym(t)$. It follows that
\[a(G)=t+a(K/N)+a(G/K)\le t+t\log{(3s\log{q})}+4(t-1)/3,\]
where the last inequality follows from Proposition \ref{prop:TransCompLength} and the paragraph above. Since $s,q\geq 2$, we quickly see that $t+t\log{(3s\log{q})}+4(t-1)/3$ is less than $t\log{s}+(s^t-1)\log{q}$ when $t\geq 4$.

We may therefore assume that $t=1,2$ or $3$. Hence, $a(G/K)\le 0,1$, or $2$, respectively. Replacing the bound $a(G/K)\le 4(t-1)/3$ with these improved bounds in the inequalities above then gives us what we need,
unless $(s,t,q)\in\{
		(2,2,4),(2,2,5),(2,2,7),(2,1,q)(4\le q\le 29),(3,1,2),(3,1,3)\}$
(note here that $(s,q)\not\in\{(2,2),(2,3)\}$, since $\GL(2,2)$ and $\GL(2,3)$ are soluble).

Since $s\in\{2,3\}$ and the relevant values of $q$ are small, we can deal with these cases by hand by examining the nonabelian simple subgroups of $\mathrm{PGL}_s(q)$ for each of the listed values of $q$ and $s$. In particular, we see that $A/S$ must be cyclic of order $1$, $2$ or $3$ in each case. It follows that $a(G)\le t+t+a(G/K)=2t+a(G/K)$. Inserting the bound $a(G/K)\le 0,1$, or $2$ when $t=1,2$, or $3$, respectively, we get what we need in each case.
	\end{proof}
	
	We are now ready to prove Theorem \ref{PrimBound}.
	\begin{proof}[Proof of Theorem \ref{PrimBound}]
		We remark first that if $(n,q)\in\{(2,3),(2,5),(4,3)\}$, then we can use \MAGMA to check that the only weakly quasiprimitive subgroups $G$ of $\GL(n,q)$ failing to satisfy $a(G)\le \log{n}+n\log{q}$ are those subgroups listed in the statement of the theorem. So we will assume throughout the proof that $(n,q)\not\in\{(2,3),(2,5),(4,3)\}$.
		
		Now, let $G$ be a counterexample to the theorem with $n$ minimal, and write $V$ for the natural $\mathbb{F}_q[G]$-module.
		By Lemma~\ref{lem:slred}, there exists a divisor $e$ of $n$ and a subnormal
		subgroup $K$ of $G$ such that each of the following holds:
		\begin{enumerate}[\upshape(1)]
			\item $a(G)\le a(K)+\log{e}$;
			\item $K$ is isomorphic to a subgroup of $\GL(n/e,q^e)$; and
			\item $K/K\cap Z(\GL(n/e,q^e))$ is isomorphic to a subdirect subgroup of a direct product $H_1\times \cdots\times H_u$, where $H_i$ satisfies one of the following:
			\begin{enumerate}[\upshape(a)]
				\item $H_i$ has shape $N.L$, where $N$ is elementary abelian of order $r_i^{2m_i}$, with $r_i$ a prime dividing
$q^e-1$, and $L$ is a subgroup of $\Sp(2m_i,r_i)$ acting naturally on $N$; or
				\item There exists a finite simple group $S$ with an absolutely irreducible projective representation $\rho$ of dimension $s_i\geq 2$ over $\mathbb{F}_{q^e}$, and a subgroup $A$ of the stabilizer of $\rho$ in $\Aut S$, such that $H_i$ is a subgroup of $A\wr \Sym(t_i)$ containing $S^t$.
			\end{enumerate}
			Moreover, setting $f(H_i):=r_i^{m_i}$ if $H_i$ is as in (a), and $f(H_i):=s_i^{t_i}$ if $H_i$ is as in (b), we have that
$\prod_{i=1}^u f(H_i)$ divides $n/e$. Finally, if $r$ is a fixed prime, then there is at most one value of $i$ in the range $1\le i\le t$ with the property that $H_i$ is as in (a) with $r_i=r$. 
		\end{enumerate}
		Now
\begin{eqnarray*}a(K) &=& a(K/K\cap Z(\GL(n/e,q^e)))+a(K\cap Z(\GL(n/e,q^e)))\\
&\le&a(K/K\cap Z(\GL(n/e,q^e)))+e\log{q},
\end{eqnarray*}
since $|Z(\GL(n/e,q^e))|=q^e-1$. Hence, by (3) above, we have 
		\begin{align*}\label{lab:Gabound}
			a(K)&\le e\log{q}+\sum_{i=1}^u a(H_i). &(8)   
		\end{align*}
		Suppose first that $u=1$. Then the result follows from Propositions \ref{prop:rmCase} and \ref{prop:stCase}, unless $H_1$ is as in (3)(a) above, and $(m_1,r_1,q^e)=(1,2,3),(1,2,5),(1,2,7)$, $(1,3,4),(2,2,3)$, or $(2,2,5)$. So assume that we are in one of these cases. Then $e=1$, $G=K$, and we can use \MAGMA to check that $a(\C_{q-1})+a(H_1)\le 5,6,6,7,10$, or $11$ respectively. Since $a(G)=a(K)\le a(\C_{q-1})+a(H_1)$ and $r_1^{m_1}$ divides $n$, this gives $a(G)\le \log{n}+n\log{q}$,
unless $r_1^{m_1}=n$ and $(m_1,r_1,q)\in\{(1,2,3),(1,2,5),(2,2,3)\}$. Since we are assuming that $(n,q)\not\in\{(2,3),(2,5),(4,3)\}$, this completes the proof of the theorem when $u=1$.
		
		Suppose now that $u>1$, and let $\mathcal{X}$ be the set of those $(m_i,r_i,q^e)$ such that there exists $1\le i\le u$ with $H_i$ as in (3)(a) above and $f(H_i)=r_i^{m_i}$. Let $\mathcal{Y}:=\{(1,2,3),(1,2,5),(1,2,7),(1,3,4),(2,2,3),(2,2,5)\}$. Since each fixed prime $r$ can be equal to $r_i$ for at most one value of $i$, the set $\mathcal{X}\cap \mathcal{Y}$ has cardinality at most $1$.
		
		Assume first that 
		$\mathcal{X}\cap\mathcal{Y}$ is not $\{(1,2,3)\}$ or $\{(2,2,3)\}$. Then by Propositions \ref{prop:rmCase} and \ref{prop:stCase}, together with (8), we have
\[ a(K)\le e\log{q}+\sum_{i=1}^u (\log{f(H_i)}+g(H_i)e\log{q}),\]
where $g(H_i)\in\{f(H_i),f(H_i)-1\}$, and $g(H_i)=f(H_i)$ if and only if $H_i$ is as in (3)(a) above and $\mathcal{X}\cap\mathcal{Y}=\{(m_i,r_i,q^e)\}$. In particular, $g(H_i)=f(H_i)-1$ for all but at most one value of $i$. Since $u>1$, we deduce that
\[a(K)\le \sum_{i=1}^u (\log{f(H_i)}+f(H_i)e\log{q}).\]
Moreover, since $f(H_i)>1$ for all $i$, we have $\sum_{i=1}^uf(H_i)\le \prod_{i=1}^tf(H_i)\le n/e$. Since $a(G)\le \log{e}+a(K)$, the theorem follows. 
		
		Finally, suppose that $\mathcal{X}\cap\mathcal{Y}$ is either $\{(1,2,3)\}$ or $\{(2,2,3)\}$. Note then that $q^e=3$, so we must have $e=1$ and $G=K$. Without loss of generality, assume that $\mathcal{X}\cap\mathcal{Y}=\{(m_1,r_1,q)\}$ (that is, $H_1$ is as in (3)(a) above with $f(H_1)=2^1$ or $2^2$). Also, set $m:=\prod_{i=2}^uf(H_i)$. Then by arguing as in the previous paragraph, we quickly see that $\log{q}+\sum_{i=2}^ua(H_i)\le\log{m}+m\log{q}$. Hence, $a(K)\le \log{q}+\sum_{i=1}^u a(H_i)\le a(H_1)+\log{m}+m\log{q}$ by (8).
Note that, since $\GL(2,3)$ is solvable and $|\mathcal{X}\cap\mathcal{Y}|=1$,
we cannot have $m=2$.

Suppose first that $(r_1,m_1,q)=(1,2,3)$.
Then $f(H_1)=2$ (so $2m$ divides $n$),  $a(H_1)\le 5$.
and $a(G)=a(K)\le 5+\log{m}+m\log{3}$.
Since $5+\log{m}+m\log{3}\le \log{2m}+2m\log{3}$ for all $m\geq 3$, we deduce that $a(G)\le \log{2m}+2m\log{3}\le \log{n}+n\log{q}$, as required.

Otherwise $\mathcal{X}\cap\mathcal{Y}=\{(2,2,3)\}$, in which case
$f(H_1)=4$, $a(H_1) \le 10$, and $a(G)=a(K)\le 10+\log{m}+m\log{3}$.
Again we find that
$10+\log{m}+m\log{3}\le \log{4m}+4m\log{3} \le  \log{n}+n\log{q}$
for $m \ge 3$, as required.
	\end{proof}

\section{Estimating the constants}\label{sec:constants}
In this final section, we attempt to estimate the constants involved in
Theorem~\ref{thm:main}.  Many of these arguments require machine computations
to check small cases. These were done in \MAGMA \cite{Magma},
but they are not difficult, and could equally easily be done in \GAP \cite{GAP}.

	\begin{lemma}\label{lem:nsbd}
		Let $G \le \GL(n,q)$ be irreducible and almost simple mod scalars, and
		suppose that $G$ does not lie in Aschbacher class $\cC_8$  (that is, it is not
		contained in the normalizer of a quasisimple classical group in its natural
		representation).
		Then $\log |G| \le \frac{5}{8} n^2 \log q$, except for the cases
		$G = A_7 \le \GL(4,2)$, $2.A_5Z \le \GL(2,9)$, and $2.A_5Z \le \GL(2,11)$,
		with $Z := Z(\GL(n,q))$, and in those three cases we have
		$\log |G| \le 0.71 n^2 \log q$.
	\end{lemma}
	\begin{proof}
		The assumptions imply that $G$ lies in one of the Aschbacher classes
		$\cC_3$, $\cC_5$ and $\cS$. The structure of the groups in classes
		$\cC_3$ (subgroups of $\Gamma L(n/e,q^e)Z$ for divisors $e>1$ of $n$) and
		$\cC_5$ (subgroups of $\GL(n,q^{1/e})Z$ for $e>1$) is well understood,
		and it it is routine to verify that the largest values of
		$\log |G|/(n^2 \log q)$ arise from the $\cC_5$-subgroups
		$\GL(2,q^{1/2})Z$ with $n=2$ and $q$ square, and these values tend to
		the supremum  $5/8$ for large $q$.
		
		By \cite{Liebeck}, Aschbacher class $\cS$-groups have order at most $q^{3n}$
		except for certain representations of $A_m$ and $S_m$ with $m \in \{n+1,n+2\}$.
		In the $A_m$ and $S_m$ cases we have $\log |G| \le 0.6 \log q$, with the
		largest value arising with $S_6 \le \GL(4,2)$.
		
		We can use the tables in \cite{BHR} to check the result for class $\cS$-groups
		when $n \le 4$, and for $n \ge 5$ Liebeck's result gives us
		$\log |G| \le 0.5 n^2 \log q$.
	\end{proof}

\subsection{The constants in the proof of Proposition~\ref{prop:primitive} and the proof of Part 3 of Theorem~\ref{thm:main}}
In this subsection we refer back to the numbered inequalities in the proof of
Proposition~\ref{prop:primitive}, and estimate the constants involved.

\paragraph{(1)}
Note first that when $e$ is prime and, in particular, when $e < 4$,
we have $d(H) \le d(H \cap K)+1$, and so we can replace the terms $C+2$ and
$2 \log \frac{n}{e} + 3$ in inequality (1) by $C+1$ and
$2 \log \frac{n}{e} + 2$, respectively.
 
With that modification, we claim that the right hand side of inequality (1) is
at most $Cn^2 \log q$ for all prime powers $q$ and all $n,e \ge 2$ with $e|n$,
provided that $C \ge 2$, which proves the result in this case.
			
To see this, note first that it is not difficult to show that the expression,
both in its original and modified form, is decreasing in $e$, so we can take
$e=4$ in general and $e=2$ with the modified expression.
In both cases it is again easy to see that we need only consider the case $q=2$.

When $e=2$ the expression becomes $n^2(C+1)/2 + 2(\log (n/2) + 1)$,
so we need to check $n^2/2 \ge 2(\log (n/2) + 1)$, which is true for $n \ge 2$.
When $e=4$, it becomes $n^2(C+2)/4 + 4\log (n/4) + 6$,
so we have to check $n^2 \ge 4\log (n/4) + 6$, which is true for $n \ge 4$.
	
\paragraph{(2)} We shall prove inequality (2) with $C_2=3$, and
$C_1 = \frac{3 \log(24)}{4\log(3)} \simeq 2.167$, where this bound is attained
with $n=2$, $q=3$, $G=\GL(2,3)$ and $H=Q_8$.
		
Furthermore, in the solvable case with $u=1$ (when $f(G_1)$ divides $n$),
we shall prove that
$d(G)\log |G| \le 2n^2 \log q$ and $d(H)\log |G| \le C n^2 \log q$ with
$C = \frac{3 \log(96)}{4\log(5)} \simeq 2.127$ (as in the statements of
Theorem~\ref{thm:main} and Proposition~\ref{prop:primitive}), where this bound is attained
when $n=2$, $q=5$ and $|G|=96$, and $H$ is a symplectic type group of order
	$16$.
			
In the non-solvable case  with $u=1$, we shall prove that
$d(H)\log |G| \le 2 n^2 \log q$.

\paragraph{(3)} We find that inequality (3) holds with $C_1 = 2$, whenever
$q \ge 8$ or $n \ge 13$, so our claimed results hold in all but
finitely many cases.  
For the exceptional cases, we use the more accurate bound on $|G_i|$ resulting
from $|G_i| \le r^{2m}.\Sp(2m,r)$ when $r$ is odd or $q \equiv 1 \bmod 4$
and $|G_i| \le r^{2m}.\GO^\pm(2m,r)$ when $r=1$ and $q \equiv 3 \bmod 4$.
Then we find that the bound holds with $C_1=2$ except when
$(f(G_i),q) = (2,3), (2,5), (3,4), (4,3)$ or $(4,5)$.
In these cases we can compute all subgroups and verify the claims directly.

\paragraph{(4)} We find that inequality (4) holds with $C_2=2$ except when
$s=t=2$.  In that case we have $S = \PSL(2,r)$
for some prime power $r$ (equal to either a root of $q$ or to $5$), and
$\PSL(2,r)^2 \le G_i\le \PGammaL(2,r) \wr \C_2$,
Its subnormal subgroup $H_i$ must either be trivial or contain
$\PSL(2,r)^2$, and we find that $d(H_i) \le 2$ for all possible $H_i$.
Now, by using the fact that $q \ge 4$ when $s=2$, we find that
the result holds with $C_2=2$.

\paragraph{(5)} For inequality (5), note that, when $S$ is an Aschbacher
$\cC_8$-group (simple classical group in its natural representation), we have
$d(H_i) \le 2$ and $d(G) \le 2$ if $u=1$, so all claims hold. So suppose not.
			
In the case $u=1$, $G$ is a nearly simple group, and 
Lemma~\ref{lem:nsbd} gives $\log|G| \le \frac{2}{3}n^2 \log q$, giving
$d(G)\log |G| \le 2n^2 \log q$, except in three specific examples.
Since $d(H) \le 3$ in general, and $d(H)$ = 2 in each of the exceptional cases,
we have $d(H)\log|G| \le 2n^2 \log q$ in all cases, which proves the claims.
	
Otherwise (when $u>1$), Lemma~\ref{lem:nsbd} gives
$\log|S| \le \frac{5}{8}f(G_i)^2 \log q$ except in three examples, and since
$G_i/S$ is isomorphic to a subgroup of $\Out S$ and
$|\Out S| \le 6\log |S|/7$ by \cite[Lemma 3.2]{Menezes}, this gives
$(d(H_i)+1)(\log(q-1) + \log |G_i|) \le
	4(\frac{5}{8}f(G_i)^2 \log q + \log(6/7))$
	(or $(d(H_i)+1)(\log(q-1) + \log |G_i|) \le
	3(0.71 f(G_i)^2 \log q + \log(6/7))$ in the exceptional case),
and it can be checked that this is less than $3 f(G_i)^2 \log q$ in all cases.

\paragraph{(6)} For inequality (6), we find that
$d(H)\log |G| \le 2n^2 \log q$ in all cases.
Using $\max(C_1,C_2) = 3$, the final inequality gives this when $u>2$, and the
penultimate inequality gives it when $u=2$, except when $f(G_1)=f(G_2)=2$.
When $u=2$ and $f(G_1)=f(G_2)=2$, we derive the result by considering the
solvable and non-solvable cases for $G_1$ and $G_2$ individually, noting that
$|G_i| \le 24$ in the solvable case.
	
\subsection{The constants in Section~\ref{sec:general} and the proof of
Part 2 of Theorem~\ref{thm:main}}\label{subsec:gencon}
	
It remains to prove Part 2 of Theorem~\ref{thm:main}.
	The following lemma, which essentially follows from \cite{Tracey18}, bounds the minimal number of elements required to generate a large subgroup in a wreath product $H\wr T$.
		\begin{lemma}\label{lem:dGbound}
			Let $H$ be a finite group, let $T\le \Sym(t)$ be a permutation group of degree $t\geq 2$, and let $G$ be a large subgroup of the wreath product $H\wr T$. Set $b:=\sqrt{2/\pi}$. Then
			\begin{enumerate}[\upshape(i)]
				\item We have $$d(G)\le \dfrac{a(H)bt\sqrt{2}}{\sqrt{\log{t}}}+d(T).$$
				\item Suppose that $T$ contains $\Alt(t)$. Then
				$$d(G)\le 2a(H)+1.$$
				\item Suppose that $H$ is a large subgroup of a wreath product $R\wr D$, where $D\in\{\Alt(d),\Sym(d)\}$ with $d\geq 3$. Then
				$$d(G)\le 2a(R)t+\dfrac{bt\sqrt{2}}{\sqrt{\log{t}}}+d(T).$$
			\end{enumerate}
		\end{lemma}
		\begin{proof}
			Write $p_1^{a_1},\hdots,p_k^{a_k}$ for the orders of the abelian chief factors of $H$, and let $E(t,p_i)$ and $E_{\sol}(t,p_i)$ be as defined in \cite[Definition 4.21]{Tracey18}, and note that $E_{\sol}(t,p_i)\le E(t,p_i)$. Write $c_{\nonab}(H)$ for the number of nonabelian chief factors of $H$. Then by \cite[Corollary 5.10(i)]{Tracey18}, we have $d(G)\le \sum_{i=1}^k a_iE(t,p_i)+c_{\nonab}(H)+d(T)$. Moreover, \cite[Corollary 4.27]{Tracey18} implies that $E(t,p_i)$ is bounded above by $bt\sqrt{2}/\sqrt{\log{t}}$ for $t\ge 1261$. One can check using the definition of $E(t,p)$ that this is also true for $2\le t\le 1260$. Since $\sum_{i=1}^ka_i+c_{\nonab}(H)\le a(H)$, part (i) follows.
			
			We now prove (ii). So assume that $\Alt(t)\le T$. Then $d(T)\le 2$. If $t\le 4$, then $T$ is soluble, so $d(G)\le \sum_{i=1}^k a_iE_{\sol}(t,p_i)+c_{\nonab}(H)+2$ by \cite[Corollary 5.10(i)]{Tracey18}. We can also compute from the definition of $E_{\sol}$ that $E_{\sol}(t,p_i)\le 2$ for $t\in\{2,3,4\}$. Since $c_{\nonab}(H)+1\le 2c_{\nonab}(H)$, and $\sum_{i=1}^ka_i+c_{\nonab}(H)\le a(H)$, the result follows.
			So assume that $t\geq 5$. Then \cite[Proposition 4.30(i) and Corollary 5.9]{Tracey18} imply that $d(G)\le \sum_{i=1}^k2a_i+c_{\nonab}(H)+2$. The result then follows as above.  
			
			Finally, we prove (iii). So assume that $H$ is a large subgroup of a wreath product $R\wr D$, where $D\in\{\Alt(d),\Sym(d)\}$ with $d\geq 3$. By \cite[Proposition 4.30(i) and Proposition 6.9 and its proof]{Tracey18}, we have $d(G)\le \sum_{i=1}^k 2a_it+c_{\nonab}(R)+1+iE(t,2)+d(T)$, where $i:=0$ if $D=\Alt(d)$, and $i:=1$ otherwise. Since $\sum_{i=1}^k 2a_it+c_{\nonab}(H)+1\le 2a(H)t$, and $E(t,2)\le bt\sqrt{2}/\sqrt{\log{t}}$ by the first paragraph of the proof, the result follows.  
	\end{proof}

\begin{proof}[Proof of Theorem~\ref{thm:main}, Part 2] 
Since Theorem~\ref{thm:main} Part 2 follows from Part 3 in the case $G$ primitive, we may assume that $G$ is imprimitive. Then as mentioned in Section~\ref{sec:general}, we may view $G$ as a large subgroup of a wreath product $H\wr T$, where $H\le \GL(s,t)$ is primitive, and $T\le \Sym(t)$ is transitive. Write $K:=G\cap H^s$ for the intersection of $G$ with the base group of $H\wr T$. Then $d(G)\log{|G|}=d(G)\log{|K|}+d(G)\log{|T|}$. As in Section~\ref{sec:general}, we will bound each of the quantities $d(G)\log{|K|}$ and $d(G)\log{|T|}$ spearately, beginning with the latter.

	
Now, when bounding $d(G) \log |T|$ in Section~\ref{sec:general}, we
identified $T$ with a subgroup of $T_1 \wr \cdots \wr T_m$,
where each $T_i$ is a primitive permutation group of degree $t_i > 1$,
and we considered first the case when none of the the groups $T_i$ with
$t_i \ge 3$ contains $\Alt(t_i)$. In that case, we have
 $|T| \le c^t$ with $c:= 2.8349$, by \cite[Corollary 1.4]{Maroti02}, and so
$\log{|T|}\le t\log{c}$.
		Assume first that $(s,q)\not\in\{(2,3),(2,5),(4,3)\}$. Then $a(H)\le\log{s}+ s\log{q}$ by Theorem \ref{PrimBound}. Also, setting $c_1:=\sqrt{3}/2$, we have $d(T)\le c_1t/\sqrt{\log{t}}$, by \cite{Tracey21Sharp}. It follows from Lemma \ref{lem:dGbound} and the above that $d(G)\log{|T|}\le \log{c}(t^2/\sqrt{\log{t}})(b\sqrt{2}\log{s}+bs\sqrt{2}\log{q}+c_1)$. Set $\beta_0:=0.6164$. When $s\geq 3$, the bound noted yields $d(G)\log{|T|}\le \beta_0 s^2t^2\log{q}$ for $t\geq 7$, and all values of $q$. When $2\le t\le 6$, then we can use \MAGMA to check that $d(T)\le i$, where $i=2$ if $t>2$, and $i=1$ if $t=2$. So we have $d(G)\log{|T|}\le \log{c}(t^2/\sqrt{\log{t}})(b\sqrt{2}\log{s}+bs\sqrt{2}\log{q})+it\log{c}$. This again yields $d(G)\log{|T|}\le \beta_0 s^2t^2\log{q}$ apart from some small values of $s$ and $q$ (and it can be easily checked by hand that the bound $d(G)\log{|T|}\le \beta_0 s^2t^2\log{q}$ still holds in these cases). A similar argument deals with the cases $s\le 2$, and $(s,q)\in\{(2,3),(2,5),(4,3)\}$.
	
	Otherwise, we chose $k$ with $1 \le k \le m$ such that $t_k$ is maximal with the property that
	$t_k \ge 3$ and $\Alt(t_k) \le T_k$. Suppose first that $k=m$. An easy argument then shows that $G$ has a block $W_{m-1}$ of dimension $r:=st_1\hdots t_{m-1}$.
	Write $H_{m-1}:=\Stab_G(W_{m-1})^{W_{m-1}}\le \GL(W_{m-1})$ for the image of the induced action of $\Stab_G(W_{m-1})$ on $W_{m-1}$.
	Then $d(G)\le 2a(H_{m-1})+1$, by Lemma \ref{lem:dGbound}(ii).
	
	Set 
	$$\beta:=10\log{(3)}/9=1.76107\hdots\text{ and }\gamma:=\beta-1/6.$$
	We claim that $d(G)\log{|T|}\le \gamma n^2\log{q}$. To see this, suppose first that $q\not\in\{3,5\}$. Then $a(H_{m-1})\le r\log{q}+4(r-1)/3$ by Corollary \ref{cor:PrimBoundcorollary}. Arguing as above, this yields $d(G)\log{|T|}\le (2r\log{q}+8r/3-1/3)t\log{t_m}= (2r\log{q}+8r/3-1/3)(n/s)\log{t_m}$. Since $n=rt_m$, $q\geq 2$, and $t_m\geq 3$, one can then quickly see that if either $s\geq 2$, $t_m\geq 10$, or $q\geq 7$, the claim follows. So assume that $s=1$, $q\le 5$, and $t_m\le 9$. Then $q>2$, since $G$ is irreducible, so $q=4$, since we are assuming that $q\not\in\{3,5\}$. It follows that $R$ is a large subgroup of $\C_3\wr T_0$ for some transitive $T_0\le \Sym(r)$. Hence, $a(R)\le r+4(r-1)/3$ by Proposition \ref{prop:TransCompLength}. We then have $d(G)\log{|T|}\le (7r/3-1/3)n\log{t_m}$, and this quickly yields $d(G)\log{|T|}\le \gamma n^2\log{4}$, as needed.

		Suppose next that $q=3$. Then $a(R)\le 19r/6-4/3$ by Corollary \ref{cor:PrimBoundcorollary}. Moreover, if $s=1$, then $R$ is a large subgroup of $\C_2\wr T_0$, where $T_0\le\Sym(r)$ is transitive. It follows from Proposition \ref{prop:TransCompLength} that $a(R)\le r+4(r-1)/3$ in this case. Replacing the bounds for $a(R)$ in the paragraph above with these new bounds for $a(R)$, and arguing in the same way, we get $d(G)\log{|T|}\le \gamma n^2\log{q}$. If $q=5$, then $a(R)\le 22r/6-4/3$ by Corollary \ref{cor:PrimBoundcorollary}, and the argument is entirely similar.
		
		We have proved that $d(G)\log{|T|}\le \gamma n^2\log{q}$ when $k=m$. Suppose next that $k<m$. We will prove that $d(G)\log{|T|}\le \beta n^2\log{q}$. To this end, note first that $G$ has a block $W_{k}$ of dimension $rt_k$, where $r:=st_1\hdots t_{k-1}$.
	Write $H_{k}:=\Stab_G(W_{k})^{W_{k}}\le \GL(W_{k})$ for the image of the induced action of $\Stab_G(W_{k})$ on $W_{k}$. It is then routine to check that $H_k$ is a large subgroup of a wreath product $R\wr T_k$, where $R$ is an irreducible subgroup of $\mathrm{GL}(r,q)$. Further, $G$ is a large subgroup of $H_k\wr \widetilde{T}$, where $\widetilde{T}$ is transitive of degree $\widetilde{t}:=t_{k+1}\hdots t_m$. Since $T_k\in\{\Alt_{t_k},\Sym(t_k)\}$, we then have
	$d(G)\le 2a(R)\widetilde{t}+b\widetilde{t}\sqrt{2}/\sqrt{\log{\widetilde{t}}}+d(\widetilde{T})$, by Lemma \ref{lem:dGbound}(iii). Also, setting $c_1:=\sqrt{3}/2$ as above, we have $d(\widetilde{T})\le c_1\widetilde{t}/\sqrt{\log{\widetilde{t}}}$, by \cite{Tracey21Sharp}. Since $|T|\le t_k^t$ by \cite[Proposition 1.9]{Tracey16}, we deduce that
	\begin{align*}\label{lab:alg}
	d(G)\log{|T|}&\le t\log{t_k}\left(2a(R)\widetilde{t}+
\frac{(b\sqrt{2}+c_1)\widetilde{t}}{\sqrt{\log{\widetilde{t}}}} \right).  &(9)  
	\end{align*}
	Suppose first that $q\not\in\{3,5\}$. Then by Corollary \ref{cor:PrimBoundcorollary}, we have $a(R)\le 4(r-1)/3+r\log{q}$. If $t\le n/2$, then since $n=rt_k\widetilde{t}$, the bound $d(G)\log{|T|}\le \beta n^2\log{q}$ follows quickly from (9). Assume, then, that $t=n$. 
    Then we can again use (9) to deduce that $d(G)\log{|T|}\le \beta n^2\log{q}$ if either $q\geq 13$ or $t_k\geq 15$. So assume further that $q\le 11$ and $t_k\le 14$. 
	Now, since $t=n$, we have that $G$ is a large subgroup of $C\wr T$ for some $C\le \C_{q-1}$. It follows that $d(G)\le n+2n/3=5n/3$, whence $d(G)\log{|T|}\le (5/3)n^2\log{t_k}$. If $t_k\le q$, then the bound $d(G)\log{|T|}\le \beta n^2\log{q}$ follows. So assume that $t_k>q$.
	Note also that $q\neq 2$, since $G$ is irreducible. (Thus, $q\geq 4$ since we are assuming that $q\neq 3$). Note that $R$ is isomorphic to a large subgroup of $\C_0\wr T_0$, for some $\C_0\le \C_{q-1}$, and some transitive $T_0\le \Sym(r)$. Thus, $a(R)\le r\omega{(q-1)}+(4/3)(r-1)$ by Proposition \ref{prop:TransCompLength}, where $\omega(q-1)$ denotes the number of prime divisors of $q-1$, counted with multiplicities.  For each choice of $t_k\le 14$ and $q\le 11$ with $q\not\in\{2,3,5\}$ and $q<t_k$, the bound $d(G)\log{|T|}\le \beta n^2\log{q}$ then follows from (9).
	
	The cases $q\in\{3,5\}$ are dealt with in an entirely similar way, using (9) and the bounds on $a(R)$ from Corollary \ref{cor:PrimBoundcorollary}.
	
We complete the proof by bounding $d(G) \log |K|$. Set $C:= \frac{3 \log(96)}{4\log(5)}$, and
suppose first that $(T,t)\neq (\Sym(3),3)$ (so that the bound $d(T)\le t/2$ holds). Then from (7) in Section~\ref{sec:general}, we see that $d(G)\log{|K|}\le (C+1/2)n^2\log{q}$. It follows from our work above that $d(G)\log|G| \le D n^2 \log q$ with
		$D \le C + \beta + 1/2=4.38806\hdots$. 
		
		Suppose finally that $(T,t)= (\Sym(3),3)$. Then we see from (7) in Section~\ref{sec:general} and our analysis in the case $k=m$ above (and from the fact that $d(T)=2t/3$ in this case) that $d(G)\log|G| \le D_1 n^2 \log q$ with
		$D_1 \le C + \gamma+ 2/3=C+\beta+1/2=4.38806\hdots$. This completes the proof.
\end{proof}

	\textsc{Derek Holt,
		Mathematics Institute,
		University of Warwick,
		Coventry CV4 7AL, UK}
	
	\emph{E-mail address}{:\;\;}\texttt{D.F.Holt@warwick.ac.uk}
	\bigskip
	
	\textsc{Gareth Tracey,
		School of Mathematics,
		University of Birmingham, Edgbaston, Birmingham, B15 2TT, UK}
	
	\emph{E-mail address}{:\;\;}\texttt{garethtracey1@gmail.com}
	
\end{document}